\documentclass[graybox]{svmult}

\usepackage{mathptmx}       
\usepackage{helvet}         
\usepackage{courier}        
\usepackage{type1cm}        
\usepackage{makeidx}         
\usepackage{graphicx}        
\usepackage{multicol}        
\usepackage[bottom]{footmisc}

\usepackage{amsmath,amssymb}

\makeindex             

\begin{document}

\title*{Modeling phase transition and metastable phases}
\author{Fran\c cois James and H\'el\`ene Mathis}
\institute{Fran\c cois James \at MAPMO, Univ. Orl\'eans
  and CNRS, UMR CNRS 7349,
  45067 Orléans Cedex 2, \email{francois.james@univ-orleans.fr}
\and H\'el\`ene Mathis \at LMJL, Univ. Nantes, 2 rue de la
Houssini\`ere - BP 92208 F-44322 Nantes Cedex3, 
\email{helene.mathis@univ-nantes.fr}}

%
%
\maketitle

\abstract{We propose a model that describes phase transition
  including meta\-stable phases present in the van der Waals Equation of State (EoS).
  We introduce a dynamical system that is able to depict the mass
  transfer between two phases, for which equilibrium states are both metastable and stable states,
including mixtures.
  The dynamical system is then used as a relaxation source 
  term in a isothermal two-phase model.
  We use a Finite volume scheme (FV) that
  treats the convective part and the source term in a fractional step
  way. Numerical results illustrate the ability of the model to capture
  phase transition and metastable states.}

\section{Introduction}
\label{sec:introduction}
Metastable liquids  are liquid states where the temperature is higher
than the ebullition temperature. Such states are very unstable and a
very small perturbation brings out a bubble of vapor inside the
liquid. Such phenomenon can appear at saturated temperature (or at
saturated pressure for metastable vapor) for instance inside a nozzle
such as fuel injector or cooling circuit of water pressurized reactor.
In the last decades considerable research has been devoted to the modeling of
two-phase flows with phase transition. 
However the exact expressions of the
transfer mass term are usually unknown (see \cite{drew83}). In
particular, to our knowledge,
there is very few literature about the transfer term able to depict
metastable states. In \cite{saurel08} and \cite{zein10} the authors
consider a 6 equation model where relaxation to equilibrium is
achieved by chemical and pressure relaxation terms whose kinetics are
considered infinitely fast.

We intend here to provide a new model able to depict phase transition and
metastable states with non-infinite relaxation speed.
It is based on the use of the van der Waals EOS,
that is well-known to depict stable and metastable states below the
critical temperature. However this EOS is not valid in the so-called
spinodal zone where the pressure is a decreasing function of the
density.
This leads to instabilities and computational failure and the pressure
has to be corrected using the Maxwell equal area rule
construction to recover a constant pressure. But such a correction
removes the metastable regions.
We propose transfer terms obtained through an optimization problem
of the Helmholtz free energy of the two-phase system.
For sake of simplicity we assume the system to be isothermal.
We obtained a dynamical system that is able to depict mass transfer
including metastable states and that dissipates the total Helmholtz free 
energy.
The equilibria of the dynamical system are both stable and metastable
states and mixture states that satisfies the pressures and
chemical potentials equalities.
This dynamical system is used as transfer term in a isothermal
two-phase model in the spirit of \cite{guillard05} and \cite{chalons09}.
We use a classical FV scheme that treats the convective and the source
terms in a splitting approach.

Section~\ref{sec:thermo} is devoted to the thermodynamics of binary
mixture and presents the major properties of the van der Waals EoS. 
Section~\ref{sec:dynam-syst-phase} is devoted to the construction of
the dynamical system based on results of the previous
Section. In particular we show that metastable states are
attractors of the dynamical system.
In Section~\ref{sec:numer-scheme-comp} we briefly present the splitting
FV scheme we use and give numerical results where metastable vapor appears.

\section{Thermodynamics and van der Waals Equation of State}
\label{sec:thermo}

In this Section we first recall the thermodynamics theory for a single
isothermal fluid and introduce the different potentials of the van der
Waals EoS, then we state the mathematical framework for the thermodynamics of
immiscible binary mixtures.

\subsection{Thermodynamics of a single phase}
\label{sec:therm-single-phase}
Consider a single fluid of mass $M>0$ occupying a volume $V>0$.
At constant temperature if the fluid is homogeneous and at rest, its
behavior is entirely described by the Helmholtz free energy function $E(M,V)$
which belongs to 
$C^2(\mathbb R_+ \times \mathbb R_+)$ and is positively homogeneous of degree 1 (PH1).
Thus, at fixed volume $V$, one can introduce the specific Helmholtz
free energy $f$ and the specific energy $e$ that are functions of the
density $\rho=M/V$ 
	\begin{equation}\label{eq:fe}
f(\rho) = E(\rho,1), \qquad \rho e(\rho) =E(\rho,1).
	\end{equation}
We introduce also the pressure $p$ and the chemical potential $\mu$
that are partial derivatives of the free energy $E$, respectively with respect to $V$ and $\rho$.
By homogeneity, one can write them as functions of $\rho$ solely:
	\begin{equation}\label{eq:pmu}
p(\rho) = -\partial_V E(\rho,1), \qquad  \mu(\rho)=\partial_M E(\rho,1).
	\end{equation}
Again thanks to the homogeneity of the energy function, one has
	\begin{equation}\label{eq:gibbs_relation}
  f(\rho) = \rho\mu(\rho)-p(\rho),\qquad  f'(\rho) = \mu(\rho).
	\end{equation}

\begin{figure}[t]
\sidecaption
\includegraphics[scale=.35]{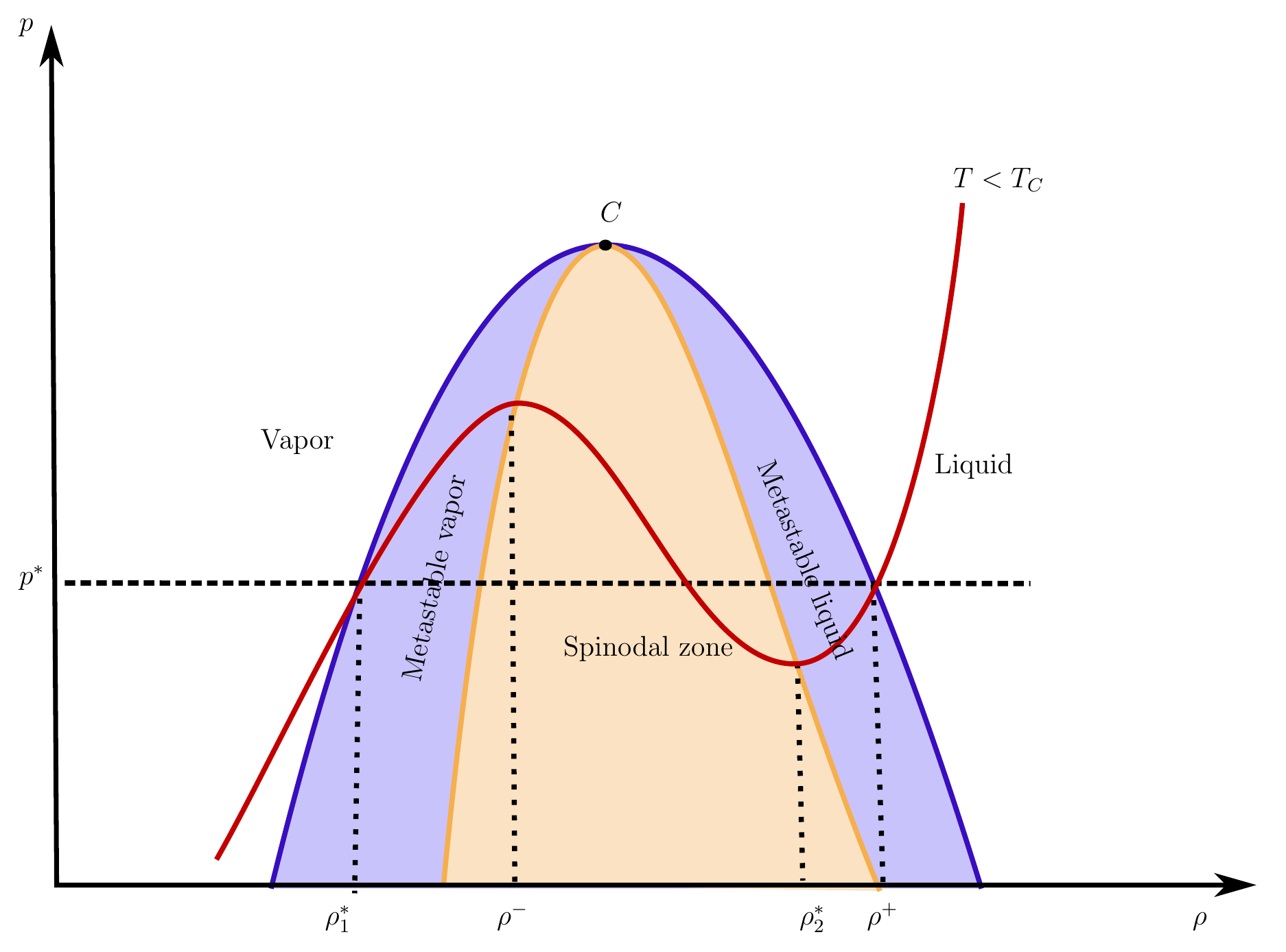}
\caption{Phase diagram for the van der Waals EoS in the $(p,\rho)$
  plan. The red curve stands for an isotherm below the critical
  temperature $T_C$, the point $C$ being the critical point. 
  The orange zone is called spinodal zone, it corresponds
  to unstable states. In that area the isotherm has to be replaced by
  an horizontal segment that coincides with the isobaric line at
  constant pressure $p^*$.}
	\label{fig:1}       
	\end{figure}

Stable pure phases are characterized by a convex energy function, which leads to a nondecreasing
pressure law. We consider a classical example of a fluid that may experience phase transitions, namely
the van der Waals monoatomic fluid. At fixed temperature $T$ its Helmholtz free energy is given by
	\begin{equation}\label{eq:E_vdW}
E(M,V) = -\dfrac{a M^2}{V} + RT\left(M\log\dfrac{M}{V-Mb}-M\right),
	\end{equation}
where $R$ stands for the perfect gas constant and $a$ and $b$ are 
positive constants, $a$ accounts for binary interactions and $b$ is the covolume. 
Below a critical temperature $T_C$ the pressure law is no longer monotone (see fig. \ref{fig:1}): in a region
called the spinodal zone, the pressure decreases with respect to the density,
thus leading to instable states. In that region the isotherm have to be
replaced by the  maxwell area rule in order to recover that phase
transition happens at constant pressure and chemical potential.
However this construction removes admissible regions where
the pressure law is still nondecreasing. Such regions are called the metastable
regions (blue regions in fig. \ref{fig:1}). 
We consider in the following the dimensionless equation of state and the
associated potentials for which $R=8/3$, $a=3$ and $b=1/3$, for which $T_C=1$.

\subsection{Equilibrium of a two-phase mixture}
\label{sec:equil-two-phase}

We consider now two immiscible phases of a same pure fluid of total
mass $M$ and volume $V$.
Each phase $i=1,2,$  is depicted by its mass $M_i\geq 0$ and its volume
$V_i\geq 0$. We assume that both phases are characterized by
the same van der Waals extensive Helmholtz free energy $E$ function of $M_i$ and
$V_i$, given by \eqref{eq:E_vdW}.
By the conservation of mass, the mass of the binary system is
$M=M_1+M_2$ and immiscibility implies $V=V_1+V_2$.

According to the second principle of thermodynamics (see
\cite{Gibbs}), for fixed mass $M$ and volume $V$
the stable equilibrium states of the system are the
solutions to the constrained optimization problem
\begin{equation*}
  \inf \{ E(M_1,V_1) + E(M_2,V_2) | \; V_1+V_2=V, \; M_1+M_2 = M\},
\end{equation*}
which  can be rewritten using~\eqref{eq:fe} in term of the specific
Helmholtz free energy at fixed density $\rho$:
	\begin{equation}\label{eq:pb_opti}
\inf \{ \alpha_1f(\rho_1) + \alpha_2f(\rho_2) | \; \alpha_1+\alpha_2=1,
	\; \alpha_1 \rho_1+\alpha_2 \rho_2 =\rho \},
	\end{equation}
where $\alpha_i=V_i/V \in [0,1]$ denotes the volume fraction and
$\rho_i=M_i/V_i$ is the density of the
phase $i=1,2$. In the sequel the fractions $\alpha_i$ are written as
functions of $\rho,\rho_1$ and $\rho_2$ such that
$\alpha_1  (\rho,\rho_1,\rho_2) =(\rho-\rho_2)/(\rho_1-\rho_2)$
and 
$\alpha_2  (\rho,\rho_1,\rho_2) = 1-\alpha_1  (\rho,\rho_1,\rho_2)$.

Note that $\alpha_1$ and $\alpha_2$ are simultaneously non zero if and
only if $\rho_1 \neq \rho_2$. In that case we shall always assume without loss of generality
that $\rho_1<\rho_2$ and $\rho\in [\rho_1,\rho_2]$.
The total Helmholtz free energy $F : \mathbb R_+^3\to \mathbb R$ of the binary system is given by
	\begin{equation}\label{eq:F}
F(\rho, \rho_1,\rho_2) = \alpha_1 (\rho, \rho_1,\rho_2)f(\rho_1) + \alpha_2 (\rho, \rho_1,\rho_2)f(\rho_2).
	\end{equation}
Depending on the saturation of the volume fractions, one can
characterize the equilibria of the optimization problem
\eqref{eq:pb_opti}.
\begin{proposition}\label{prop:eq}
  \begin{enumerate}
\item Pure states: if $\alpha_1=0$ (resp. $\alpha_2=0$) then only the phase $2$
    (resp. 1) is stable.
\item Mixture: if $\alpha_1 \alpha_2 \neq 0$, then the equilibrium state
is characterized by one of the following  equivalent properties
    \begin{enumerate}
\item equality of the chemical potentials and the pressures
    \begin{equation}
      \label{eq:eq_p_mu}
      \mu(\rho_1) = \mu(\rho_2)=\mu^*,\quad
      p(\rho_1) =p(\rho_2)=p^*,
    \end{equation}
\item Maxwell area rule on the chemical potential
    \begin{equation}
      \label{eq:maxwell_mu}
      \int_0^1 \mu(\rho_2 + t(\rho_1-\rho_2))dt = \mu(\rho_1) = \mu(\rho_2)=\mu^*,
    \end{equation}
\item the difference of the energies reads
    \begin{equation}
      \label{eq:e2-e1}
      f(\rho_2) - f(\rho_1) = \mu (\rho_1)(\rho_2-\rho_1) =
      \mu(\rho_2) (\rho_2-\rho_1).
    \end{equation}
    \end{enumerate}
  \end{enumerate}
The densities such that \eqref{eq:eq_p_mu}, \eqref{eq:maxwell_mu} or \eqref{eq:e2-e1}
hold  are denoted $\rho_1^*$ and $\rho_2^*$, see fig. \ref{fig:1}.
\end{proposition}
The most important consequence of this result is that in the metastable zones there are two possible equilibrium states
corresponding to a pure metastable state and a stable mixture state. Hence the EoS at equilibrium is not single-valued. The difference
between stable and metastable states lies in their dynamical behaviour with respect to perturbations, see \cite{Landau}.

\section{Dynamical system and phase transition}
\label{sec:dynam-syst-phase}

We turn now to the study of dynamical stability of equilibrium states. First we address the homogenous case, 
introducing a dynamical system for which the equilibria 
are both stable and metastable states as well as states in the spinodal
area such that~\eqref{eq:eq_p_mu}-\eqref{eq:maxwell_mu} are satisfied.
Next the dynamical system is plugged as a relaxation 
source terms in a isothermal two-fluid model. Some properties of the
full model are given: hyperbolicity, existence of a energy function
that decreases in time.
 
\subsection{Dynamical system}
\label{sec:dynamical-system}
Assuming that $\rho$, $\rho_1$ and $\rho_2$ are only time-dependent, 
we introduce the following dynamical system, which derives from the optimality conditions of
Proposition \ref{prop:eq}:
	\begin{eqnarray}\label{eq:syst_dyn}
\dot{\rho} &= &0,\nonumber\\
\dot\rho_1 &= & -(\rho-\rho_1)(\rho-\rho_2)\left(
      \rho_2(\mu(\rho_2)-\mu(\rho_1)) + p(\rho_1) -p(\rho_2)\right),\\
\dot\rho_2 &= & (\rho-\rho_1)(\rho-\rho_2)\left(
      \rho_1(\mu(\rho_1)-\mu(\rho_2)) - p(\rho_1) +p(\rho_2)\right).\nonumber
	\end{eqnarray}
Straightforward computions show that the total Helmholtz free energy $F$
defined by \eqref{eq:F} decreases in time along the solutions
of this system.
We focus now on the equilibria which can be reached by the model
(under the assumption $\rho_1<\rho_2$).
\begin{theorem}
  \label{th:equi}
  The equilibria of the system~\eqref{eq:syst_dyn} are
  \begin{enumerate}
  \item the monophasic states such that $\alpha_1=1$ (resp. $\alpha_1=0$) that is
    $\rho=\rho_1=\overline\rho$ with any
    $\rho_2\neq \overline\rho$ (resp. $\rho=\rho_2=\overline\rho$ with any
    $\rho_1\neq \overline\rho$).
    In that case, if $\rho=\rho_i$, $i=1$ or $2$ such that
    \begin{enumerate}
    \item $\rho\not\in[\rho^-, \rho^+]$, then the equilibrium is an 
      attractor and corresponds to
      monophasic and metastable states,
    \item $\rho \in[\rho^-, \rho^+]$, then the equilibrium is a repeller and corresponds to
      states belonging to the spinodal zone (which is non admissible),
    \end{enumerate}
  \item  the unique state such that $0<\alpha_1<1$ and the
    relations~\eqref{eq:eq_p_mu}-\eqref{eq:maxwell_mu} are satisfied. 
  \end{enumerate}
\end{theorem}
A remarkable feature of this system is that a perturbation of a pure metastable state involving the other phase leads
to a mixture equilibrium state, corresponding to the definition of metastable state \cite{Landau}.

\subsection{The isothermal model}
\label{sec:isothermal-model}

The previous dynamical system \eqref{eq:syst_dyn} is now coupled with
a modified version of the isothermal two-phase model proposed in
\cite{chalons09} (see also \cite{guillard05}). 
The model admits
a mixture pressure $\alpha_1 p(\rho_1) + \alpha_2 p(\rho_2)$ and one
velocity $u$ for both phases.
It reads
\begin{equation}
  \label{eq:model}
  \begin{aligned}
    \partial_t \rho + \partial_x (\rho u) &= \dfrac 1 \varepsilon\dot
    \rho =0,\\
    \partial_t \rho_i + \partial_x(\rho_i u )&= \dfrac 1 \varepsilon\dot{\rho_i},\qquad i=1,2\\
    \partial_t (\rho u) + \partial_x (\rho u^2 + \alpha_1 p(\rho_1) +
    \alpha_2 p (\rho_2)) &= 0,
  \end{aligned}
\end{equation}
where the source terms are given by the dynamical system
\eqref{eq:syst_dyn} and account for mass and mechanical transfer.
The parameter $\varepsilon>0$ is a relaxation parameter that represents
the relaxation time to reach thermodynamical equilibrium. In order to
capture metastable states, we will consider $0<\varepsilon<1$ in computations.

The convective part of the model \eqref{eq:model} is hyperbolic with
the eigenvalues
\begin{equation}
  \label{eq:eig}
  \lambda_1=u-c, \quad, \lambda_2=\lambda_3=u, \quad \lambda_4=u+c,
\end{equation}
where the speed of sound is $c=\sqrt{\dfrac 1 \rho
  (\alpha_1 \rho_1 p'(\rho_1) + \alpha_2 \rho_2 p'(\rho_2))}$.
\begin{proposition}
  The function $\mathcal E (\rho, \rho_1, \rho_2, u) = \dfrac{\rho u^2}{2} +
  \alpha_1 f(\rho_1) + \alpha_2 f(\rho_2)$, satisfies the following equation
  \begin{equation}
    \label{eq:comp_bal_law}
    \partial_t (\mathcal E) + \partial_x (u(\mathcal E + \alpha_1
    p(\rho_1) + \alpha_2 p(\rho_2)) = (\partial_{\rho_1} F)\dot{\rho_1} +
    (\partial_{\rho_2} F)\dot{\rho_2} \leq 0.
  \end{equation}
\end{proposition}
Note that $\mathcal E$ is not an entropy of the system since $f$ is
a non-convex function of the density.

\section{Numerical illustration}
\label{sec:numer-scheme-comp}

We present here numerical results that assess the ability of the model to
capture phase transition including metastable states. We use a standard Finite Volume method to
approximate the Cauchy problem 
	\begin{equation}\label{eq:cauchy}
\partial_t W + \partial_xF(W) = S(W),\qquad W(0,x) = W_0(x), \, x \in \mathbb R,
	\end{equation}
where $ W =(\rho, \rho_1,\rho_2, \rho u)^T$,
$F(W) =(\rho u, \rho_1 u,\rho_2 u, \rho u^2 + \alpha_1p(\rho_1) +
\alpha_2 p(\rho_2))^T$,
and 
$S(W) = (0,\dfrac 1 \varepsilon\dot\rho_1, \dfrac 1
\varepsilon\dot \rho_2, 0)^T$.
We use a fractional step approach. We denote $\Delta t$ the time step and
$\Delta x$ the length of the cell $(x_{i-1/2}, x_{i+1/2})$ on the
regular 1D-mesh. Let $W^n$ be the Finite Volume approximation at time
$t^n = n\Delta t$, $n\in \mathbb N$. 
The first step corresponds to the approximation of the convective part
which provides the solution $W^{n,-}$ at time $t^{n,-}$. It is treated by a classical 
Rusanov scheme. The second step is the approximation of the source terms (relaxation),
at this stage we merely use an explicit Euler method.

\begin{figure}[ht]
  \centering
  \includegraphics[height=2.7cm]{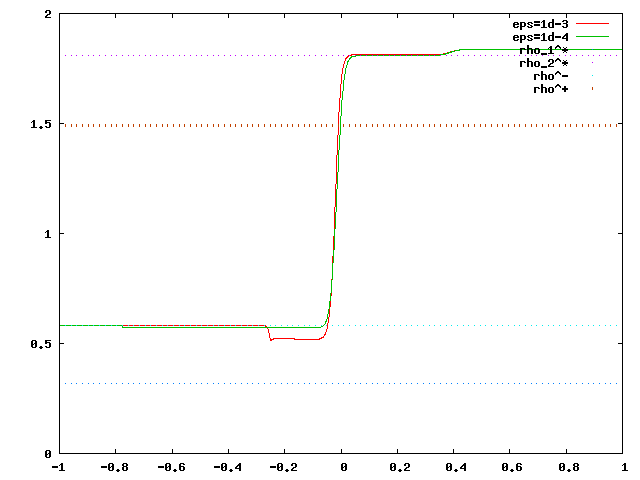}
  \includegraphics[height=2.7cm]{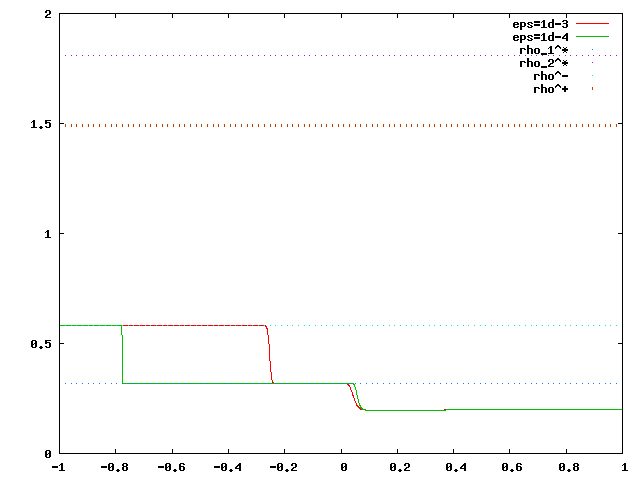}
  \includegraphics[height=2.7cm]{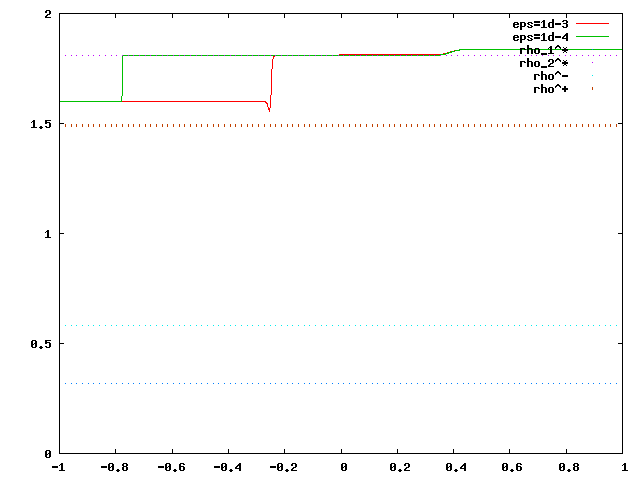}
  \includegraphics[height=2.7cm]{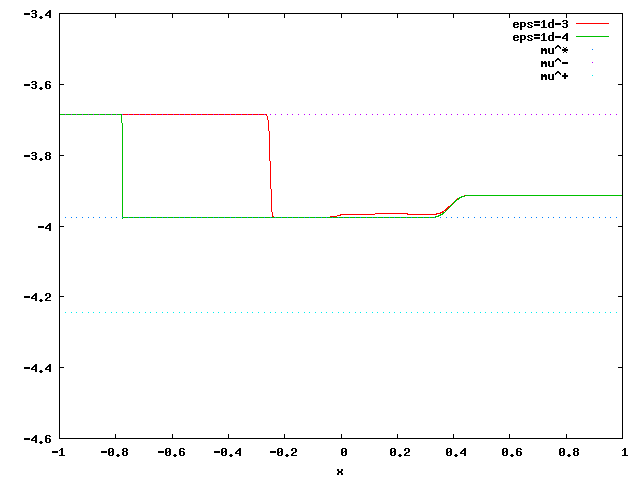}
  \includegraphics[height=2.7cm]{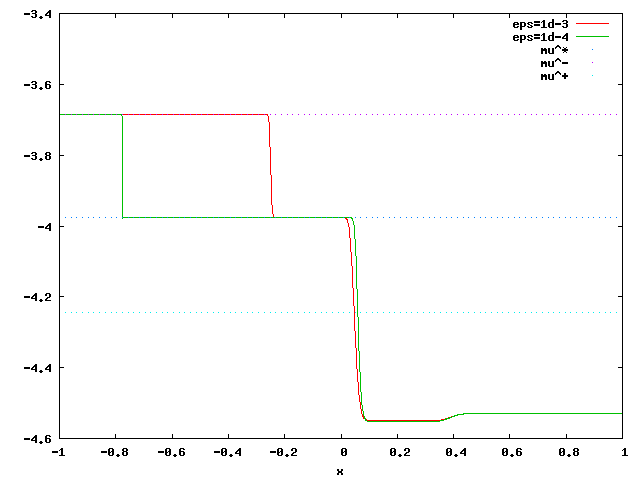}
  \includegraphics[height=2.7cm]{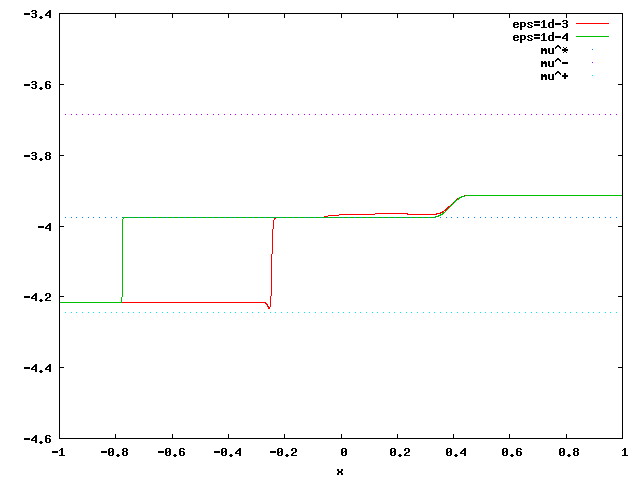}
  \includegraphics[height=2.7cm]{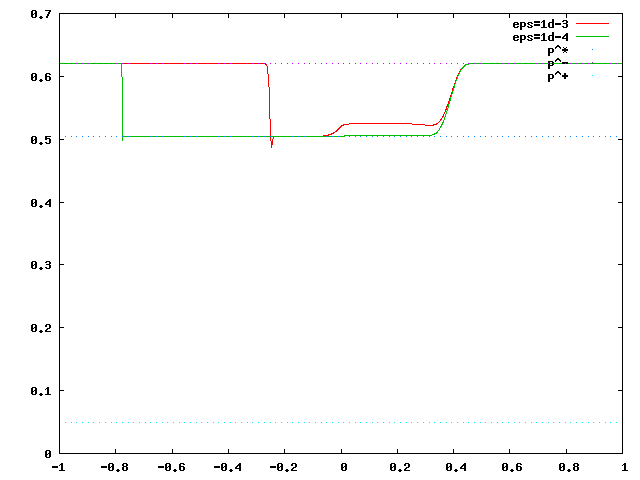}
  \includegraphics[height=2.7cm]{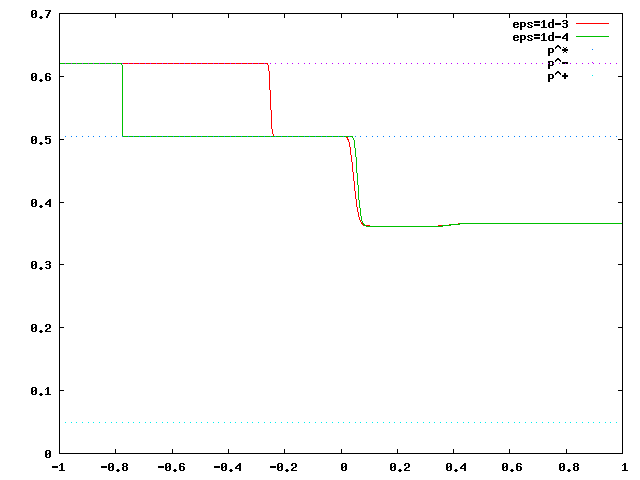}
  \includegraphics[height=2.7cm]{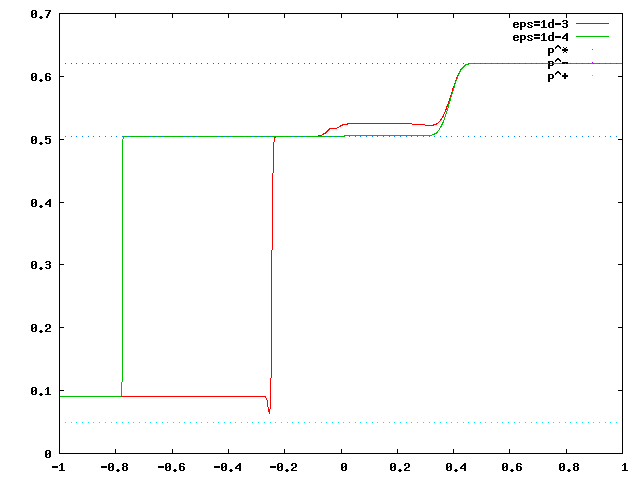}
  \includegraphics[height=2.7cm]{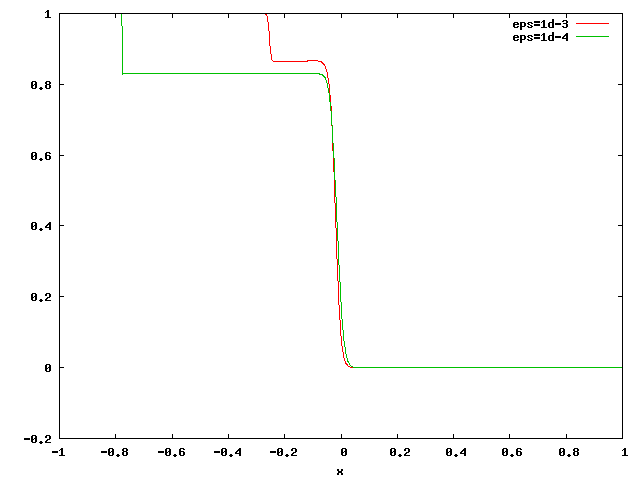}
  \includegraphics[height=2.7cm]{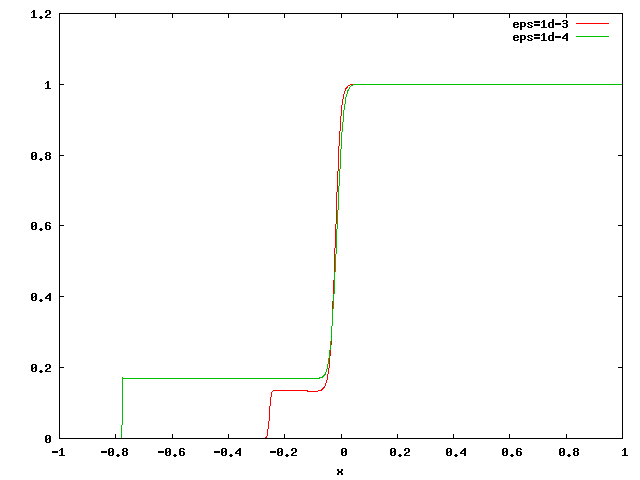}
  \includegraphics[height=2.7cm]{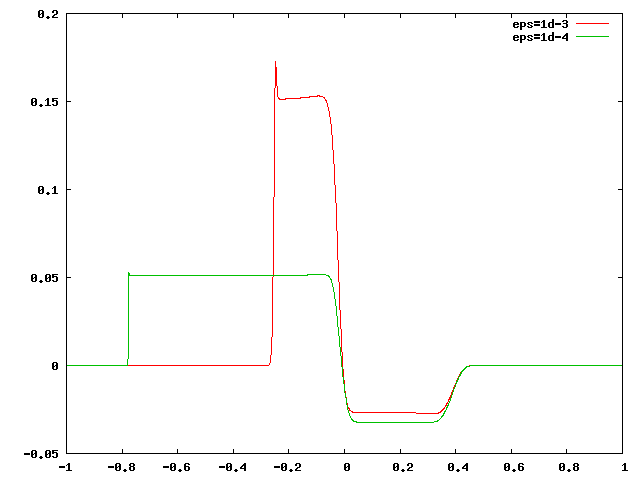}
  \caption{ First line, densities: $\rho$, $\rho_1$,
  $\rho_2$, second line, chemical potentials:
  $\mu$, $\mu_1$, $\mu_2$, third line, pressures: $p$, $p_1$, $p_2$, last line
  fractions $\alpha_1$, $\alpha_2$ and the velocity $u$.}
  \label{fig:test2_annexe2}
\end{figure}

We consider the van der Waals equation at constant temperature
$T=0.85$. The extrema of the isotherm curve are
$\rho^-=0.581079$ and $\rho^+=1.488804$. The Maxwell construction
on the chemical potential defines the densities
$\rho_1^* = 0.319729$ and $\rho_2^*=1.807140$ 
such that $\mu(\rho_1^*)=\mu(\rho_2^*)=3.977178$
and $p(\rho_1^*)=p(\rho_2^*)=0.504492$. 
If the Riemann problem consists in an initial constant pressure and
constant chemical potential state, the numerical scheme preserves
this state exactly as it is expected.
Another test case consists in an initial constant pressure state which is
subjected to a disequilibrium in chemical potential.
The initial data are $\rho_L=\rho_{1,L}=\rho^-$, $\rho_{2,L}=1.6$,
$\rho_R=\rho_{2,R}=1.837840$, $\rho_{1,R}=0.2$ and $u_L=u_R=0$.
The discontinuity is applied at $x=0$ in the domain $[-1,1]$. The mesh
contains with 2000 cells and the time of computation is $t=0.2$.
Note that $\rho_{2,L}$ belongs to the metastable liquid region and
$\rho_{2,R}$ belongs to the pure liquid region such that
$p(\rho_{2,R}) = p(\rho^-)=p(\rho_{1,L})$ and $\rho_{1,R}$
belongs to the pure gaseous region.
Fig. \ref{fig:test2_annexe2} presents the results for $\varepsilon=10^{-3}$ and
$\varepsilon=10^{-4}$. 
The main feature to notice here is that the relaxation approximation introduces
a mixture zone on both sides of the interface, which remains stable. Within this zone, there are variations of the velocity,
which remains compressive ($u>0$ for $x<0$, $u<0$ for $x>0$).

\section{Conclusion and prospects}
\label{sec:conclusion}

The first tests with this model show that it is able to cope with phase transitions with metastable
states using a van der Waals EoS. 
Due to the complexity of the source term, we propose as a first step
an explicit treatment of the relaxation term. We aim at providing a
semi-implicit scheme in the spirit of \cite{GHS04}.
Moreover this model is a toy one since it is isothermal. We attend to
add the temperature dependance to obtain a fully heat, mass and
mechanical transfer model in order to compare our results to the one
of \cite{saurel08} and \cite{zein10}.

\noindent\textbf{Acknowledgement}
The second author is supported by the project~ANR-12-IS01-0004-01
GEONUM.

\bibliographystyle{spmpsci}
\bibliography{fvca7-james-mathis-template}
\end{document}